\input amstex
\documentstyle{amsppt}
\mag=\magstep1
\pageheight{20.5cm}
\vcorrection{-0.2cm}
\NoBlackBoxes

\document
\def\SO{\text{\rm SO}} \def\Sp{\text{\rm Sp}}  
\def\U{\text{\rm U}} \def\PP{\Bbb P} \def\R{\Bbb R} \def\C{\Bbb C}  
\def\H{\Bbb H} \def\Z{\Cal Z} \def\CP{\Cal P}  
\def\FS{\text{\rm FS}} \def\sp{\text{\rm sp}} \def\so{\text{\rm so}} \def\id{\text{\rm id}} 
\def\CY{\text{\rm CY}} \def\can{\text{\rm can}} \def\std{\text{\rm std}} \def\Lie{\text{\rm Lie}} 
\def\ov{\overline} \def\wt{\widetilde} \def\Ric{\text{\rm Ric}} \def\tr{\text{\rm tr}} 
\def\p{\partial} \def\Rm{\text{\rm Rm}}

\topmatter
\title Moving Frames on the Twistor Space of Self-Dual Positive Einstein 4-Manifolds \endtitle
\author Ryoichi Kobayashi and Kensuke Onda \endauthor
\affil Graduate School  of Mathematics, Nagoya University \endaffil
\abstract{The twistor space $\Z$ of self-dual positive Einstein manifolds 
naturally admits two 1-parameter families of Riemannian metrics, one is the family of 
canonical deformation metrics and the other is the family introduced by B. Chow 
and D. Yang in [C-Y]. The purpose of this paper is to compare these two families. 
In particular we compare the Ricci tensor and the behavior under the Ricci flow of 
these families. As an application, we propose a new proof to the fact that 
a locally irreducible self-dual positive Einstein 4-manifold is isometric to either $S^4$ 
with a standard metric or $\PP^2(\C)$ with a Fubini-Study metric.} 
\endabstract
\rightheadtext{Moving Frames on the Twistor Space}
\leftheadtext{}
\endtopmatter

\beginsection 0. Introduction
\par

Let $\R^4$ be the oriented Euclidean 4-space which is identified with $\H$ 
and we consider the right multiplication of the group $\Sp(1)$ of unit quaternions. 
This defines a subgroup $\Sp(1)_-$ in $\SO(4)$ which is isomorphic to $\Sp(1)$. The centralizer 
of this subgroup in $\SO(4)$ is again isomorphic to $\Sp(1)$ identified with 
the left multiplication of unit quaternions, which we denote by $\Sp(1)_+$. 
We thus have the group homomorphism 
$\Sp(1)\times \Sp(1) \rightarrow \SO(4)$ defined by the left/right multiplication of unit 
quaternions on $\H$. Its kernel is $\pm(\id,\id) \cong\Bbb Z_2$ and therefore we have 
the decomposition of the Lie algebra $\so(4)=\sp(1)_+\oplus\sp(1)_-$. 
Here, $\sp(1)_{\pm}$ are copies of $\sp(1)$ and $\sp(1)_+$ (resp. $\sp(1)_-$) 
corresponds to the left (resp. right) multiplication of unit quaternions. 
In terms of the identification $\so(4)=\Lambda^2\R^4$, the component $\sp(1)_+$ 
(resp. $\sp(1)_-$) corresponds to the self-dual (resp. anti-self-dual) 2-forms. 

A Riemannian 4-manifold $(M,g_M)$ is said to be self-dual positive Einstein 
if its self-dual part $W_-$ of the Weyl tensor vanishes 
and Einstein condition is satisfied with positive scalar curvature\footnote{\,\,The twistor 
space of a self-dual positive Einstein 4-manifold is a positive K\"ahler-Einstein manifold 
and hence simply connected. This implies that any self-dual positive Einstein 4-manifold 
is simply connected.}. 
The twistor space of a self-dual Riemannian 4-manifold $(M,g)$ is defined as follows. 
Let $\CP$ be the holonomy reduction of the bundle of oriented orthonormal frames 
of $(M,g)$. For instance, if $M=(S^4,g_{\std})$,  $\CP$ is a principal $SO(4)$-bundle 
over $S^4$ and if $M=(\PP^2,g_{\FS})$, $\CP$ is a principal $S(U(1)\times U(2))$-bundle. 
The holonomy group of a self-dual positive Einstein manifold $(M,g)$ is a subgroup of 
$\SO(4)$ and we write it as $H$. 
The fiber of $\CP$ over $m\in M$ consists of the $H$-rotation of an 
oriented orthonormal frame of the tangent space $M_m$. 
The twistor space $\Z$ of $(M,g)$ is defined as
$$\Z=\CP/H\cap \U(2)$$
and the twistor fibration is the $\PP^1$-bundle 
$$\pi:\Z \rightarrow M\,\,.$$
Note that the subgroup of $H$ which induces a holomorphic transformation 
w.r.to the orthogonal complex structure represented by a point of $\PP_m^1$ 
is identified with $H\cap \U(2)$ and therefore $\PP_m^1$ is identified with 
$$H/H\cap U(2)=\SO(4)/\U(2)=\Sp(1)_+\Sp(1)_-/\U(2)=\Sp(1)_-/\Sp(1)_-\cap\U(2)\,\,.$$
As $H\cap\U(2)$ is precisely the subgroup of $H$ consisting of 
all rotations in $H$ which are holomorphic w.r.to 
an orthogonal complex structure of $M_m$ (say, that induced from the right $j$-multiplication), 
the $\PP^1$-fiber of the twistor fibration is identified with the set of all orthogonal complex 
structures of $M_m$. 
The meaning $\U(2)$ is explained as follows.  Take $j$ as an orthogonal complex structure. 
The subgroup of $\SO(4)$ which centralizes the right multiplication of $j$ 
is the product (in $\SO(4)$) of $\Sp(1)_+$ and the subgroup of $\Sp(1)_-$ 
which centralizes $j$, i.e., the $S^1$-subgroup of $\Sp(1)_-$ 
consisting of $x+yj$, $x^2+y^2=1$. The product group in $\SO(4)$ is isomorphic to $\U(2)$. 
The twistor space $\Z$ of a self-dual positive Einstein 4-manifold $(M,g_M)$ 
is naturally equipped with two kinds of 1-parameter families of Riemannian metrics 
in the following way : 
\medskip

(1) The canonical deformation metrics (see, for instance, [Be]). This is defined as follows. 
The Levi-Civita connection of $(M,g_M)$ canonically defines 
a horizontal distribution on the twistor fibration $\pi:\Z \rightarrow M$. 
This allows one to define
$$g_{\lambda}^{\can}:=\lambda^2 g_{\FS}+g_M$$
where $g_{\FS}$ is Fubini-Study metric of constant curvature $4$ along the $\PP^1$-fiber 
of the twistor fibration, $g_M$ is the base self-dual positive Einstein metric on $M$ 
normalized so that its scalar curvature $48$ and $\lambda>0$ is a 
partial scaling parameter\footnote{\,\,In this paper we normalize the Fubini-Study metric 
on $\PP^1$ so that its curvature is $4$ and the self-dual positive Einstein metric on 
$M^4$ so that its scalar curvature is $48$.}. If $Tz\Z=V\oplus H$ is the decomposition into 
the vertical and horizontal subspaces, then we have $g_{\lambda}(V,H)=0$, $g_{\lambda}^{\can}|_V=\lambda^2 g_1|_V$ 
and $g_{\lambda}^{\can}|_H=g_1^{\can}|_H$.  
\medskip

(2) The Chow-Yang metrics ([C-Y]). This is very much involved compared to 
the canonical deformation metrics. 
This family was introduced by B. Chow and D. Yang in [C-Y] for the twistor space 
of positive quaternion K\"ahler manifolds. 
Roughly speaking the Chow-Yang metrics are defined as follows. 
Let $\PP^1_m$ be the $\PP^1$-fiber of the twistor fibration $\pi:\Z\rightarrow M$ 
over $m\in M$. Along $\PP^1_m$ we locally associate a smooth family of ``rotating" 
oriented orthonormal frames where the rotation is chosen from $H$ 
modulo $H\cap \U(2)$. These local associations define local sections 
$s : \Z \rightarrow \CP$ of the principal $H\cap U(2)$-bundle 
$p : \CP \rightarrow \Z$. 
Now the Chow-Yang metric is defined by requiring the dual orthonormal coframes 
of $M$ together with the image in the twistor space of the $\sp(1)_-$-part of 
the connection form (of the Levi-Civita connection of $(M,g)$) determined 
by the ``rotation" of the oriented orthonormal frames along the fiber of $\CP$ 
to be the system of orthonormal coframe on $\Z$.  
Here, the contribution from the $\sp(1)_-$-part of the connection form corresponds to 
the tangent space of the point $z$ of the fiber $\PP^1_m$ in $\Z$ which represents 
the infinitesimal deformation of the orthogonal complex structure at $z \in\PP_m^1$. 
For instance, as an orthogonal complex structure is represented by $ai+bj+ck$ 
with $a^2+b^2+c^2=1$, if the orthogonal complex structure represented by $z$ 
is the right $j$ multiplication ($j$ is defined by $a=c=0$ and $b=1$), 
then the infinitesimal deformation of the orthogonal complex structure 
is represented by $\alpha_1i+\alpha_3k$. Therefore the contribution from the 
$\sp(1)_-$-part of the connection form is $\alpha_1$ and $\alpha_3$ in this case. 
If we introduce a partial scaling 
parameter $\lambda>0$ in the connection form part, Chow-Yang metric looks like
$$g_{\lambda}^{\CY}=\lambda^2(\alpha_1^2+\alpha_3^2)+\sum_{i=0}^3X_i^2$$
where the $X_i$-part represents the ``rotating" (here, ``rotation" is considered 
modulo $H\cap \U(2)$) oriented orthonormal frame along the $\PP^1$-fiber, 
the pair $\{\alpha_1,\alpha_3\}$ (say) 
represents the connection form part representing the ``rotation modulo (orthogonal to) 
$H\cap \U(2)$-rotation" of the $X_i$-part in the $\PP^1$-fiber 
and $\lambda>0$ is a partial scaling parameter. 
\medskip

In this paper, we are concerned with the comparison of the above introduced 
two 1-parameter families of Riemannian metrics on the twistor space $\Z$ 
of self-dual positive Einstein 4-manifolds $(M,g)$. In Section 1, we compare 
these two families from topological view point and produce a semi-global Riemannian 
invariant which distinguishes these two families. Then we compute the Ricci tensor 
of metrics in these two families. The computation shows that the Ricci tensor 
of the Chow-Yang (resp. the canonical deformation) metrics is again the Chow-Yang 
(resp. the canonical deformation) type. 
In Section 2, we compare these two families from the Ricci flow view point. 

A convention. We always normalize the scaling of the self-dual positive 
Einstein 4-manifold under consideration so that its scalar curvature is $48$ (if $M=S^4$ 
with the standard metric $g_{\std}$, then its sectional curvature is identically $4$, and 
if $M=\PP^2(\C)$ with the Fubini-Study metric, then its holomorphic sectional curvature is 
identically $8$) and normalize the Fubini-Study metric of $\PP^1$ so that its Gaussian 
curvature is $4$. 

Then our main results are concerned with compact self-dual Einstein 4-manifolds 
with positive scalar curvature : 

\proclaim{Theorem 0.1} A canonical deformation metric and a Chow-Yang metric on the 
twistor space $\Z$ coincide if and only if the partial scaling parameter $\lambda$ satisfies 
the condition $\lambda=1$. In this case $g_1^{\can}=g_1^{\CY}$ is a K\"ahler-Einstein 
metric on $\Z$. 
\endproclaim

\proclaim{Theorem 0.2} (1) If $\lambda^2<3$, the Ricci tensor of the canonical deformation 
metric $g_{\lambda}^{\can}$ is positive and given by the formula
$$\Ric(g_{\lambda}^{\can})=4(1+\lambda^4)g_{\FS}+4(3-\lambda^2)g_M
=4(3-\lambda^2)\,g^{\can}_{\sqrt{\frac{1+\lambda^4}{3-\lambda^2}}}\,\,.$$
In particular $g_{\lambda}^{\can}$ is K\"ahler-Einstein if and only if $\lambda^2=1$ 
and Einstein non-K\"ahler if and only if $\lambda^2=\frac12$. 
\medskip

(2) The Ricci tensor of the Chow-Yang metric 
$g_{\lambda}^{\CY}$ is given by
$$\Ric(g_{\lambda}^{\CY})
=2\lambda^{-2}(1+3\lambda^2)\,g^{\CY}_{\sqrt{\frac{2\lambda^2(1+\lambda^2)}{1+3\lambda^2}}}\,\,.
$$
In particular, $g_{\lambda}^{\CY}$ is Einstein if and only if $\lambda^2=1$ 
and in this case $g_1^{\CY}$ is K\"ahler-Einstein. 
\endproclaim

\proclaim{Theorem 0.3} (1) Consider a 2-parameter family 
$$\Cal F^{\can}=\{\rho g_{\lambda}^{\can}\}_{\rho>0,\sqrt{3}>\lambda>0}$$
of canonical deformation metrics. 
The explicit form of the Ricci flow with initial metric in $\Cal F^{\can}$ is given by
$$\left\{\aligned
& \frac{d\lambda^2}{dt}=-\frac{8}{\rho}(\lambda^2-1)(2\lambda^2-1)\\
& \frac{d\rho}{dt}=-8(3-\lambda^2)\,\,.\endaligned\right.
$$
The trajectory is given by the equation
$$\rho=(\text{\rm const.})\,\frac{|\lambda^2-1|^2}{|2\lambda^2-1|^{\frac52}}\,\,.$$

(2) Consider a 2-parameter family 
$$\Cal F^{\CY}=\{\rho g_{\lambda}^{\CY}\}_{\rho>0,\lambda>0}$$
of the Chow-Yang metrics. 
The explicit form of the Ricci flow with initial metric in $\Cal F^{\CY}$ is given by
$$
\left\{\aligned
& \frac{d}{dt}(\rho(t)\lambda(t)^2)=-8(1+\lambda^2(t))\,\,,\\
& \frac{d}{dt}\rho(t)=-4(\lambda(t)^{-2}+3)\,\,.\endaligned\right.
$$
The trajectory is given by the equation
$$\rho=(\text{\rm const.})\,\frac{\lambda^2}{|1-\lambda^2|^4}\,\,.$$
The family $\Cal F^{\CY}$ is interpreted as a ``Ricci flow unstable cell" centered 
at a K\"ahler-Einstein metric in the sense that each solution of this family is an ancient solution 
with the K\"ahler-Einstein metric as its asymptotic soliton.
\endproclaim

The motivation of Theorem 0.3 is the gradient flow interpretation of the Ricci flow 
proposed by Perelman [P]. 
Here, a Ricci flow solution is called an ancient solution if the existence time is 
$(-\infty,T)$ ($T$ being finite) (see [Ha] and [P]). 
In the trajectory of the Ricci flow on the family $\Cal F^{\CY}$ 
(identified with the $(\lambda^2,\rho)$-plane), if $\lambda^2>1$ increases, 
then $\rho$ decreases. Moreover, in the Ricci flow in the family $\Cal F^{\CY}$ 
the scaling parameter $\rho(t)$ decreases along the flow. 
Therefore the partial scaling parameter $\lambda^2(t)$ (if $\lambda^2>1$) in 
$\Cal F^{\CY}$ increases along the flow. On the other hand, if $\lambda^2>1$, 
$\lambda^2$ decreases along the Ricci flow in the family $\Cal F^{\can}$. 

\proclaim{Theorem 0.4} A locally irreducible self-dual positive Einstein 4-manifold 
is isometric to either $(S^4,g_{\std})$ or $(\PP^2(\C),g_{\FS})$. 
\endproclaim

Here $(S^4,g_{\std})$ is $S^4$ with the standard metric of constant sectional 
curvature $4$ (scalar curvature $48$) and $(\PP^2(\C),g_{\FS})$ is $\PP^2(\C)$ 
with the Fubini-Study metric of constant holomorphic sectional curvature $8$ 
(scalar curvature $48$). 
\medskip

Theorem 0.4 is a classical result in 4-dimensional geometry 
(see Hitchin [Hi], Friedrich-Kurke [F-K] and [S2]). 
A new proof proposed in this paper is based on 
the fact that the set of the family of the Chow-Yang metrics for $\lambda>1$ 
is foliated by the trajectories of the Ricci flow ancient solutions whose asymptotic soliton 
is the trajectory of the shrinking homothetical family of the K\"ahler-Einstein metrics on $\Z$. 
The advantage of our proof lies in the fact that we can apply the 
same strategy to the positive quaternion K\"ahler case of all dimensions $\geq 8$ 
in a uniform way. In [K-O], we will extend our methods to positive quaternion K\"ahler case 
and construct a ``Ricci flow unstable cell" centered at the K\"ahler-Einstein metric 
on the twistor space of positive quaternion K\"ahler manifolds and apply this Ricci flow 
unstable cell to the study of the structure of positive quaternion K\"ahler manifolds. 
In particular we will show that any locally irreducible positive quaternion K\"ahler manifold 
is isometric to one of the Wolf spaces (this is an affirmative answer to the LeBrun-Salamon 
conjecture ([Leb], [L-S])). 

\beginsection 1. Moving Frames on the Twistor Space

Let $(M,g)$ be a compact self-dual positive Einstein 4-manifold. As introduced in the 
introduction, let $\CP$ be the holonomy reduction of the bundle of oriented orthonormal 
frames of $(M,g)$. Then the holonomy group $H(\subset \SO(4))$ acts from the right by 
rotation of a given oriented orthonormal frame. 
We define the twistor space $\Z$ of $(M,g)$ by
$$\Z=\CP/H\cap \U(2)$$
and write $\pi:\Z \rightarrow M$ for the twistor fibration. The fibration $\CP\rightarrow \Z$ 
is a principal $H\cap \U(2)$-bundle and each fiber over a point of 
$M$ is a copy of $H/H\cap \U(2)=\SO(4)/\U(2)=\Sp(1)_-/\Sp(1)_-\cap\U(2)=\PP^1$, i.e., 
the fibration $\CP\rightarrow\Z$ restricted on the $\PP^1$-fiber of the twistor fibration 
is identified with the upper or the lower row of the diagram
$$\CD
U(2) @>{}>> SO(4) @>{}>> \PP^1\\
@V{}VV       @V{}VV            @VV{\text{\rm identity}}V\\
\Sp(1)_-\cap H @>{}>> \Sp(1)_-      @>{}>>    \PP^1
\endCD$$
where the lower row is just the Hopf fibration $S^1 \rightarrow S^3 \rightarrow \PP^1$. 
The twistor space $\Z$ fibers over $M$ with fiber $\PP^1$ and the fiber $\PP_m^1$ 
over $m\in M$ is identified with the set of all orthogonal complex structures on the tangent 
spaces $M_m$. 
The Levi-Civita connection of $(M,g)$ canonically defines 
a horizontal distribution on the twistor fibration. This allows us to define a 
1-parameter family of Riemannian metrics $g_{\lambda}^{\can}$ which is the so called 
family of canonical deformation metrics 
(the sum $g_{\lambda}^{\can}=\lambda^2 g_{\FS}+g_M$ 
of the base metric $g_M$ on $M$ 
and the scaled fiber Fubini-Study metric $\lambda^2 g_{\FS}$). 

In [C-Y], B. Chow and D. Yang introduced another 1-parameter family of metrics 
(called the family of Chow-Yang metrics) $g_{\lambda}^{\CY}$ on 
the twistor space $\Z$ of positive quaternion K\"ahler manifolds by using the moving frame 
technique. 

We proceed to the description of the Chow-Yang metrics. First we consider 
topologically. 
Let $(M,g)$ be a self-dual positive Einstein 4-manifold, $\Z$ its twistor space and 
$\pi:\Z\rightarrow M$ the twistor fibration. Let $m\in M$, $U$ an open neighborhood of 
$m\in M$ and consider a local oriented orthonormal frame on $U$. 
We consider any local frame field on $U$ which defines a section $U \rightarrow \CP|_U$ 
of the $H$-principal bundle $\CP\rightarrow M$ (the holonomy reduction of 
the bundle of oriented orthonormal frames). We extend the section in the direction 
of the fiber $\CP_m$ over $m$ of $\CP\rightarrow M$ by rotating the given 
oriented orthonormal frame at $m$ by elements of $H$. 
This defines a ``local" section $s:\Z|_U \rightarrow \CP|_U$ of the principal 
$H\cap\U(2)$-bundle $\CP \rightarrow \Z$, where ``local" means 
that the section is defined locally along the $\PP^1$-fiber of the twistor fibration
\footnote{\,\,This is the unique way to extend the local orthonormal frame field on $U(\subset M)$ 
to the inverse image of $U$ under the twistor fibration $\pi:\Z \rightarrow M$ in the way compatible 
with the geometric meaning (i.e., the rotation by the left $H$-action) of $\CP$. }. 
This section corresponds to the ``rotation" of the given oriented orthonormal frame by 
elements of $H$ modulo $\U(2)$. Because there exists no global section 
of the $S^1$-bundle $H\rightarrow H/H\cap\U(2)\cong\PP^1$, 
the section is defined only locally along the $\PP^1$-fiber of the twistor fibration. 
Thus, given local oriented orthonormal frame field on $U$, we defines a ``local" section 
$s:\Z|_U\rightarrow \CP|_U$ and therefore we get the system of 1-forms on $\Z$
$$\{X^0,X^1,X^2,X^3\}$$
which is the pull-back of the unitary coframes defined globally on $\CP$ and 
the system of 1-forms on $\Z$
$$\{\alpha_1,\alpha_2,\alpha_3\}$$
which is considered modulo $H\cap\U(2)$-action and is the pull-back of 
the Levi-Civita connection form of $(M,g)$. 
For instance, if $z\in \PP_m^1$ (the $\PP^1$-fiber of the twistor fibration over $m \in M$) 
corresponds to the orthogonal complex structure $j$ represented by 
$(0,1,0)$, the infinitesimal deformation of the orthogonal 
complex structure at $j$ is given by 
$\alpha_1i+\alpha_3k$. The Chow-Yang metrics are defined by
$$
g_{\lambda}^{\CY}=\lambda^2(\alpha_1^2+\alpha_3^2)+\sum_{i=0}^3X_i^2\
$$
by declaring that the system 
$$\{\lambda \alpha_1,\lambda \alpha_3,X_0,X_1,X_2,X_3\}$$
being an orthonormal coframe ($\lambda>0$ being a partial scaling parameter). 
We are now ready to compare the family of Chow-Yang metrics and canonical 
deformation metrics. To state the result, we normalize the base space $(M,g)$ 
so that the canonical deformation metric for $\lambda=1$, i.e., $g_1^{\can}=g_{FS}+g_M$ 
is K\"ahler-Einstein. 

\proclaim{Theorem 1.1 (Theorem 0.1)} Let $(M,g)$ be a self-dual positive Einstein manifold. 
The canonical deformation metric $g_{\lambda}^{\can}=\lambda^2 g_{\FS}+g_M$ 
and the Chow-Yang metric 
$g_{\lambda'}^{\CY}=\lambda^2(\alpha_1^2+\alpha_3^2)+\sum_{i=0}^3X_i^2$ 
coincide if and only if $\lambda=\lambda'=1$. 
In particular, the Chow-Yang metric $g_{\lambda'}^{\CY}$ for $\lambda'\not=1$ 
do not belong to the family of canonical deformation metrics.
\endproclaim
\noindent
{\it Proof.} It is well-known that the canonical deformation metric is K\"ahler-Einstein 
for a unique suitable partial scaling. Our normalization is that $g_1^{\can}$ is 
K\"ahler-Einstein. In this case the parallel translation along curves in the $\PP^1$-fiber 
preserves the orthogonal complex structure of the twistor space and therefore 
the above mentioned ``rotation" along the $\PP^1$-fiber must belong to $\U(2)$. 
In [C-Y], Chow and Yang proved that $g_1^{\CY}$ is K\"ahler-Einstein in the case of 
positive quaternion K\"ahler manifolds. We will prove that $g_1^{\CY}$ is 
K\"ahler-Einstein in the case that $(M,g)$ is a self-dual positive Einstein manifold 
below by moving frame computation. This means that 
a canonical deformation metric $g_{\lambda}^{\can}$ and a Chow-Yang metric 
$g_{\lambda'}^{\CY}$ coincide if $\lambda=\lambda'=1$. 

Suppose next that $\lambda,\lambda'\not=1$. 
For $g_{\lambda}^{\can}$, there exists an oriented orthonormal frame 
field for the horizontal subspaces defined globally along a $\PP^1$-fiber 
of the twistor fibration (namely, the constant horizontal frame along the $\PP^1$-fiber). 
We show that such a global object does not exist for the Chow-Yang 
metric $g_{\lambda'}^{\CY}$ for $\lambda'\not=1$. 
To see this we fix a $\PP^1$-fiber of the twistor fibration. 
We note that for any value of $\lambda'>0$, any $\PP^1$-fiber is totally geodesic 
w.r.to the metric $g_{\lambda'}^{\CY}$. 
Therefore the parallel translation along curves in the $\PP^1$-fiber preserves 
tangent spaces of the $\PP^1$-fiber and therefore preserves the horizontal subspaces. 
On the other hand, as was shown in [C-Y] by Chow and Yang for positive quaternion 
K\"ahler manifolds and will be shown below by moving frame computation for 
self-dual positive Einstein 4-manifolds, the Chow-Yang metric $g_{\lambda'}^{\CY}$ 
is never K\"ahler if $\lambda'\not=1$. 
Therefore the holonomy along closed curves in the $\PP^1$-fiber 
is not contained in $\U(2)$. Therefore the set of all parallel translations along curves 
in the $\PP^1$-fiber of a given oriented horizontal orthonormal frame 
must be identical to $H/H\cap \U(2)\cong \PP^1$. 
This implies that there exists no smooth horizontal oriented orthonormal frame field 
defined globally along the $\PP^1$-fiber. 
Indeed, the existence of such a global object would correspond to a global section 
of the principal $H\cap\U(2)$-bundle 
$H\rightarrow H/H\cap\U(2)\cong\PP^1$ 
which is either a copy of the Hopf fibration $S^1\rightarrow S^3 \rightarrow S^2$ or 
$\U(2)\rightarrow \SO(4)\rightarrow S^2$. However, such a global section 
does not exist. 
We have thus proved that $g_{\lambda'}^{\CY}$ is never a canonical deformation metric, 
because the set of all parallel translations along curves in the $\PP^1$-fiber 
of a given horizontal oriented orthonormal frame is a Riemannian 
invariant and this set has different topological structures for two metrics $g_{\lambda}^{\can}$ 
and $g_{\lambda'}^{\CY}$ (here, $\lambda,\lambda'\not=1$). Indeed, the set has 
a global horizontal section along the $\PP^1$-fiber for the canonical deformation metric 
$g_{\lambda}^{\can}$, while this set does not admit such a global object for the 
Chow-Yang metric $g_{\lambda'}^{\CY}$. 
\qed

\bigskip

Next we proceed to the description of the Chow-Yang metrics and compute its Ricci 
tensor. 
The right multiplication of $i$ and $j$ on $\H$ is given by the matrics
$$\pmatrix 0&-1&0&0\\ 1&0&0&0\\ 0&0&0&1\\ 0&0&-1&0\endpmatrix \qquad \text{\rm and} 
\qquad \pmatrix 0&0&-1&0\\ 0&0&0&-1\\ 1&0&0&0\\ 0&1&0&0\endpmatrix\,\,.$$
Therefore the Lie algebra of $\Sp(1)_-$ (the $\Sp(1)$ subgroup of $\SO(4)$ 
stemming from the right multiplication of unit quaternions) is given by
$$
\pmatrix 0&-a_1&-a_2&-a_3\\ 
a_1&0&a_3&-a_2\\ 
a_2&-a_3&0&a_1\\ 
a_3&a_2&-a_1&0\endpmatrix\,\,.
$$
Similarly, the Lie algebra of $\Sp(1)_+$ (the centralizer of the right multiplication 
of unit quaternions) is computed as
$$
\pmatrix 0&-A_1&-A_2&-A_3\\ A_1&0&-A_3&A_2\\ A_2&A_3&0&-A_1\\ 
A_3&-A_2&A_1&0\endpmatrix\,\,.
$$
Therefore the Lie algebra of $\SO(4)$ is expressed as
$$
\pmatrix 0&-A_1-a_1&-A_2-a_2&-A_3-a_3\\ A_1+a_1&0&-A_3+a_3&A_2-a_2\\ 
A_2+a_2&A_3-a_3&0&-A_1+a_1\\ A_3+a_3&-A_2+a_2&A_1-a_1&0\endpmatrix\,\,.
$$
Let $(M,g)$ be a self-dual 4-manifold and $\CP$ the principal $H$-bundle 
over $M$ obtained by the holonomy reduction of the 
bundle of oriented orthonormal frames of $(M,g)$. 
Any local oriented orthonormal frame field $(e_A)_{A=0}^3$ defined on an 
open set $U\subset M$ defines a section $e:U \rightarrow \CP$ 
and therefore the system $\{X^A\}_{A=0}^3$ of the 
dual coframes is globally defined on $\CP$. The Levi-Civita connection of $(M,g)$ 
defines a unique $\Lie(H)$-valued 1-form which is called the connection form. 
The connection form is characterized by the following properties (1) and (2) :
\medskip

(1) Along the fiber over $m\in M$ of the fibration $p:\CP \rightarrow M$, 
the connection form is just the infinitesimal version of the ``rotation" 
given by the right action of $H$ on an oriented orthonormal frame of $M_m$ 
(this consists of the contribution from 
the left action of elements of $\Sp(1)_+$ and the right action of elements of 
$\Sp(1)_-\cap H$ on $M_m\cong\R^4$). 
\medskip

(2) If $X$ is a tangent vector and $\{e_A\}_{A=0}^3$ is a local oriented orthonormal 
frame field on an open set $U\subset M$, then
$$\nabla_Xe_B=e_A\Gamma^A_B(X)$$
where $\Gamma=(\Gamma^A_B)_{A,B=0,1,2,3}$ represents the $\sp(1)_+$-valued 
connection 1-form corresponding to the Levi-Civita connection of $(M,g)$. 
\medskip

It follows that the connection form $\Gamma$ is written as
$$
\Gamma=(\Gamma^A_B)=\pmatrix 
0&-\Gamma_1-\alpha_1&-\Gamma_2-\alpha_2&-\Gamma_3-\alpha_3\\ 
\Gamma_1+\alpha_1&0&-\Gamma_3+\alpha_3&\Gamma_2-\alpha_2\\ 
\Gamma_2+\alpha_2&\Gamma_3-\alpha_3&0&-\Gamma_1+\alpha_3\\ 
\Gamma_3+\alpha_3&-\Gamma_2+\alpha_2&\Gamma_1-\alpha_1&0
\endpmatrix\,\,.
$$
The 1-forms $X^A$ ($A=0,1,2,3$) and $\Gamma$ satisfy the following first and second 
structure equations : 
$$
\split
& dX^A+\Gamma^A_BX^B=0\,\,.\\
& d\Gamma^A_B+\Gamma^A_C\wedge\Gamma^C_B=\Omega^A_B
\endsplit
$$
where $\Omega=(\Omega^A_B)$ is the skew symmetric matrix of 2-forms on $\CP$ 
which is identified with the curvature tensor of $(M,g)$ as follows. For each point 
$e=(e_A)\in\CP$ over $m\in M$ we have
$$R(X,Y)e_B=e_A\Omega^A_B(X,Y)$$
for all $X,Y\in M_m$. The sectional curvature for the 2-plane spanned by $\{e_A,e_B\}$ 
is given by
$$K(e_A,e_B)=g(R(e_A,e_B)e_B,e_A)=\Omega^A_B(e_A,e_B)\,\,.$$

Let $(e_A)\in\CP$ be an oriented orthonormal frame of $M_m$. 
This canonically defines an identification 
$$\R^4=\H \ni (x^0+ix^1+jx^2+kx^3) \leftrightarrow (x^0+jx^2,x^1+jx^3) \in \C^2\,\,.$$
The right multiplication of $j$ on $M_m\cong \C^2$ induces the canonical almost 
complex structure on $\C^2$. The infinitesimal deformation of the orthogonal 
complex structures at $j$ is given by $\alpha_1i+\alpha_3k$ 
where $\{\alpha_1,\alpha_2,\alpha_3\}$ is the connection form of the Levi-Civita 
connection. In this situation, we pick a point $z$ in the $\PP^1$-fiber 
over $m\in M_m$ corresponding to the right $j$ multiplication. 
The canonical almost complex structure of $\Z$ at $z$ 
is defined by determining the basis of the space of $(1,0)$-forms in the following way : 
$$
\split
& \zeta^0=\alpha_1+i\alpha_3\,\,,\\
& Z^1=X^0+iX^2\,\,,\\
& Z^2=X^1+iX^3\,\,.
\endsplit
$$
Recall that the Chow-Yang metric on the twistor space $\Z$ of $(M,g)$ 
is defined by declaring that 
$$\{\alpha_1,\alpha_3\}$$
(representing the ``rotation modulo 
$H\cap \U(2)$" of oriented orthonormal frames along the $\PP^1$-fiber of the twistor fibration) 
together with 
$$\{X_0,X_1,X_2,X_3\}$$
(the rotating oriented orthonormal frame of the 
horizontal subspace, where the rotation means the $H$-rotation modulo $H\cap U(2)$-action) 
form an orthonormal frame. 
From here on we lower the indices and write $X_0,X_1,X_2,X_3$ instead of 
$X^0,X^1,X^2,X^3$ (and same for $Z$'s). 
To get a moving frame expression of the first structure equation 
of the Chow-Yang metric $g_1^{\CY}$, we begin with the first structure equation for the metric 
$g_M$. It follows from the above expression 
of the connection form $\Gamma=(\Gamma^{A}_{B})=\cdots$ that the first structure equation 
$$dX+\Gamma\wedge X=0$$
on $\CP$ is written as
$$\left\{\aligned
& dZ_1+\ov Z_2\wedge\zeta_0+(\Gamma_0+i(\Gamma_2+\alpha_2))\wedge Z_1
+(-\Gamma_1+i\Gamma_3)\wedge Z_2=0\\
& dZ_2-\ov Z_1\wedge\zeta^0+(\Gamma_1+i\Gamma_3)\wedge Z_1
+(\Gamma_0-i(\Gamma_2-\alpha_2))\wedge Z_2=0
\endaligned\right.$$
which in matrix form is expressed as
$$\split
d\pmatrix Z_1\\ Z_2\endpmatrix&=
-\pmatrix \Gamma_0+i\Gamma_2 & -\Gamma_1+i\Gamma_3\\ 
\Gamma_1+i\Gamma_3 & \Gamma_0-i\Gamma_2\endpmatrix
\pmatrix Z_1\\ Z_2\endpmatrix\\
& \quad -\pmatrix i\alpha_2 & 0\\ 0 & i\alpha_2\endpmatrix
\pmatrix Z_1\\ Z_2\endpmatrix
-\pmatrix 0 & -(\alpha_1+i\alpha_3)\\ \alpha_1+i\alpha_3 & 0\endpmatrix
\pmatrix \ov Z_1\\ \ov Z_2\endpmatrix\,\,.
\endsplit\tag1$$
To get the first structure equation of the Chow-Yang metric on $\Z$, we need 
to compute $d\zeta_0$. For the computation of $d\zeta_0$, 
we use the second structure equation
$$d\Gamma+\Gamma\wedge\Gamma=\Omega\,\,.$$
The curvature form can be expressed as
$$
\Omega=\pmatrix
0&\Omega^0_1&\Omega^0_2&\Omega^0_3\\ 
\Omega^1_0&0&\Omega^1_2&\Omega^1_3\\ 
\Omega^2_0&\Omega^2_1&0&\Omega^2_3\\ 
\Omega^3_0&\Omega^3_1&\Omega^3_2&0
\endpmatrix$$
where $\Omega^{\mu}_{\nu}=-\Omega^{\nu}_{\mu}$. Let $(\mu,\eta,\nu)$ be any cyclic 
permutation of $(1,2,3)$. Then the above expression of the connection form implies
$$\split
\Omega^{\mu}_0&=d(\Gamma_{\mu}+\alpha_{\mu})
+(\Gamma_{\mu}+\alpha_{\mu})\wedge\Gamma_0
+(\Gamma_{\eta}-\alpha_{\eta})\wedge(\Gamma_{\nu}-\alpha_{\nu})\\
&\quad +\Gamma_0\wedge(\Gamma_{\mu}+\alpha_{\mu})
+(-\Gamma_{\nu}+\alpha_{\nu})\wedge(\Gamma_{\eta}+\alpha_{\eta})\\
&=d\Gamma_{\mu}+\Gamma_{\mu}\wedge\Gamma_0
+\Gamma_{\eta}\wedge\Gamma_{\nu}+
\Gamma_{0}\wedge\Gamma_{\mu}
-\Gamma_{\nu}\wedge\Gamma_{\eta}\\
&\quad +d\alpha_{\mu}
-2\alpha_{\eta}\wedge\alpha_{\nu}\\
&=d\Gamma_{\mu}+2\Gamma_{\eta}\wedge\Gamma_{\nu}+d\alpha_{\mu}
-2\alpha_{\eta}\wedge\alpha_{\nu}\,\,\\
\Omega^{\eta}_{\nu}&=d(-\Gamma_{\mu}+\alpha_{\mu})
+(\Gamma_{\eta}+\alpha_{\eta})\wedge(-\Gamma_{\nu}-\alpha_{\nu})
+(-\Gamma_{\mu}+\alpha_{\mu})\wedge\Gamma_0\\
&+(\Gamma_{\nu}-\alpha_{\nu})\wedge(\Gamma_{\eta}-\alpha_{\eta})
+\Gamma_0\wedge(-\Gamma_{\mu}+\alpha_{\mu})\\
&=-d\Gamma_{\mu}-\Gamma_{\eta}\wedge\Gamma_{\nu}
-\Gamma_{\mu}\wedge\Gamma_0
+\Gamma_{\nu}\wedge\Gamma_{\eta}
-\Gamma_0\wedge\Gamma_{\mu}\\
&\quad +d\alpha_{\mu}
-2\alpha_{\eta}\wedge\alpha_{\nu}\\
&=-d\Gamma_{\mu}-2\Gamma_{\eta}\wedge\Gamma_{\nu}
+d\alpha_{\mu}-2\alpha_{\eta}\wedge\alpha_{\nu}\,\,.
\endsplit$$
Therefore we have
$$\Omega^{\mu}_0+\Omega^{\eta}_{\nu}
=2d\alpha_{\mu}-4\alpha_{\eta}\wedge\alpha_{\nu}\,\,.\tag2$$
To compute further, we need some representation theory of the curvature tensor 
of Riemannian 4-manifolds ([A-H-S]). 
We first normalize that the standard metric on the 4-dimensional sphere $(S^4,g_{\std})$ 
has constant sectional curvature $4$ (and therefore the scalar curvature $48$). 
We consider the curvature operator of a Riemannian 4-manifold $(M,g)$ as an element 
of an endomorphism of the bundle of 2-forms which decomposes according to the 
Lie algebra decomposition $\so(4)=\so(3)+\so(3)$ into self-dual and anti-self-dual 
2-forms : 
$$\Lambda^2(M)=\Lambda^2_+(M)\oplus\Lambda^2_-(M)\,\,.$$
According to the above decomposition of $\Lambda^2(M)$, 
the curvature operator of $(M,g)$ decomposes as
$$
\Cal R=\frac{s}{48}\Cal R_0+\Cal R' \tag3
$$
where $s$ is the scalar curvature of $(M,g)$, $\Cal R_0=\pmatrix 4\id&0\\ 0&4\id\endpmatrix$ 
is the curvature operator of $(S^4,g_{\std})$ of scalar curvature $48$ and the remaining part 
$\Cal R'$ is expressed in block matrix form as
$$\Cal R'=\pmatrix W_+&B\\ B&W_-\endpmatrix$$
where $W_+$ (resp. $W_-$) is the self-dual (resp. anti-self-dual) part of the 
Weyl tensor and $B$ is the traceless Ricci tensor. The self-dual (resp. anti-self-dual) part 
of the curvature operator consists of $\pmatrix W_++\frac{s}{12}\id&B\endpmatrix$ 
(resp. $\pmatrix B&W_-+\frac{s}{12}\id\endpmatrix$). 
In particukar, If $(M,g)$ is a self-dual positive Einstein 4-manifold, the $\Cal R'$ part 
has a matrix form
$$
\Cal R'=\pmatrix W_+&0\\0&0\endpmatrix
$$
with respect to the above decomposition of $\Lambda^2(M)$. 
Therefore $\Cal R'$ has the same form as the curvature operator of a self-dual 
Einstein manifold of zero scalar curvature (i.e., a complex K3 surface with a Ricci-flat 
K\"ahler metric orientation reversed). This part may be called the ``hyper-K\"ahler part" 
of the curvature operator. Note that the hyper-K\"ahler part has no contribution to the 
computation of the Ricci curvature of a self-dual positive Einstein manifold 
(for a quaternion K\"ahler generalization of $\dim \geq 8$ of the above story, see [S1]). 

The curvature form $\wt\Omega$ of $(S^4,g_{\std})$ of sectional curvature $4$ 
is expressed as
$$
\wt\Omega^{\kappa}_{\lambda}
=4X_{\kappa}\wedge X_{\lambda}\qquad (\kappa,\lambda=0,1,2,3)\,\,.\tag4
$$
Let $\wt\alpha_{\mu}$ ($\mu=1,2,3$) denote the $\alpha_{\mu}$'s for $(S^4,g_{\std})$ 
of curvature $4$. From (2) and (4) we have for any cyclic permutation $(\mu,\eta,\nu)$ 
of $(1,2,3)$ : 
$$
d\wt\alpha_{\mu}-2\wt\alpha_{\eta}\wedge\wt\alpha_{\nu}
=2(X_{\mu}\wedge X_0+X_{\eta}\wedge X_{\nu})\,\,.\tag5
$$
We write $\Omega'$ for the part of the curvature form corresponding to $\Cal R'$. 
Then it follows from (3) and (5) (or directly from (4)) that
$$
\split
& \quad d\alpha_{\mu}-2\alpha_{\eta}\wedge\alpha_{\nu}
=\frac12(\Omega^{\mu}_0+\Omega^{\eta}_{\nu})\\
& = \frac12\frac{s}{48}(\wt\Omega^{\mu}_0+\wt\Omega^{\eta}_{\nu})
+\frac12({{\Omega}'}^{\mu}_0+{{\Omega}'}^{\eta}_{\nu})\\
& = \frac{s}{48}(d\wt\alpha_{\mu}-2\wt\alpha_{\eta}\wedge\wt\alpha_{\nu})\\
& = \frac{2s}{48}(X_{\mu}\wedge X_0+X_{\eta}\wedge X_{\nu})
\endsplit
$$
where the third equality is a direct consequence from the fact that the $\Omega'$-part 
does not involve the $\alpha$-part. Indeed, the $\Omega'$-part (in the positive Einstein case) 
consists of $W_+$ in the above block form which is contained in the $\sp(1)_+$-part of the 
curvature form. On the other hand, the $\alpha$-part stems from the 
``rotation modulo $H\cap \U(2)$" of oriented orthonormal frames along the $\PP^1$-fiber 
of the twistor fibration. 
Since the $\sp(1)_+$-part (stemming from the left multiplication of unit quaternions) acts 
holomorphically w.r.to the orthogonal complex structures (introduced by the right multiplication 
of unit imaginary quaternions), this does not contribute to the above mentioned ``rotation" 
because the ``rotation" induced by the $W_+$-part is contained in the action of $H\cap \U(2)$. 
This means that the $W_+$ part has no contribution to the $\alpha$-part. 

From the above identity we have
$$\split
d\zeta_0&=d(\alpha_1+i\alpha_3)\\
&=2\alpha_2\wedge\alpha_3+(d\alpha_1-2\alpha_2\wedge\alpha_3)+2i\alpha_1\wedge\alpha_2
+i(d\alpha_3-2\alpha_1\wedge\alpha_2)\\
&=-2i\alpha_2\wedge\zeta_0+\frac{s}{48}(Z_2\wedge Z_1-Z_1\wedge Z_2)\,\,.
\endsplit\tag6$$
Combining (1) and (6) we have the first structure equation for the Hermitian 
metric $g_1^{\CY}$ on $\Z$ with respect to the basis of $(1,0)$-forms 
(we take the normalization $s=48$ into account) : 
$$
d\pmatrix \zeta_0\\ Z_1\\ Z_2\endpmatrix=-\pmatrix 2i\alpha_2&-Z_2& Z_1\\
\ov Z_2&i\Gamma_2+i\alpha_2&-\Gamma_1+i\Gamma_3\\
-\ov Z_1&\Gamma_1+i\Gamma_3&-i\Gamma_2+i\alpha_2\endpmatrix
\wedge \pmatrix \zeta_0\\ Z_1\\ Z_2\endpmatrix\,\,.\tag7$$
First, the right hand side contains no $(0,2)$-forms and this implies that the almost complex 
structure defined by specifying by the space $\{\zeta_0,Z_1,Z_2\}$ of $(1,0)$-forms 
(i.e., the orthogonal complex structure) is integrable. 
Second, the connection matrix is certainly skew-Hermitian and this means that the Hermitian 
metric $g_1^{\CY}$ on $\Z$ is K\"ahler with respect to the orthogonal complex structure. 

Moreover, the direct computation shows 
that the curvature form of $g_1^{\CY}$ is expressed as
$$
\pmatrix 2\zeta_0\wedge\ov\zeta_0+Z_1\wedge\ov Z_1 & \zeta_0\wedge \ov Z_1 & 
\zeta_0\wedge \ov Z_2\\
+Z_2\wedge \ov Z_2 & {}   & {} \\
{}& {}& {}\\
Z_1\wedge \ov\zeta_0& i\Omega^2_0 & 
-\frac12\{\Omega^1_0+\Omega^3_2-i(\Omega^3_0+\Omega^2_1)\}\\
{}& -\ov Z_2\wedge Z_2+\zeta_0\wedge\ov{\zeta}_0&+\ov Z_2\wedge Z_1\\
{}&{}&{}\\
Z_2\wedge \ov\zeta_0&\frac12\{\Omega^1_0+\Omega^3_2+i(\Omega^3_0+\Omega^2_1)\}
& i\Omega^3_1\\
{}&+\ov Z_1\wedge Z_2&-\ov Z_1\wedge Z_1+\zeta_0\wedge\ov\zeta_0
\endpmatrix\,\,.
$$
Therefore, the Ricci form in the complex notation is
$$\Ric(\Omega)=\tr(\Omega)=4(\zeta_0\wedge\ov\zeta_0+Z_1\wedge\ov Z_1
+Z_2\wedge\ov Z_2)\,\,.$$
This implies that $g_1^{\CY}$ is a K\"ahler-Einstein metric on $\Z$ ([C-Y]). 

On the other hand, it is well-known that the canonical deformation metric $g_1^{\can}$ 
is K\"ahler-Einstein (see, for instance, [Be]). Therefore we have
$$g_1^{\CY}=g_1^{\can}=\text{\rm K\"ahler-Einstein metric on $\Z$}\,\,.$$

Next we proceed to the study of the moving frame of the metric
$$g_{\lambda}^{\CY}=\lambda^2(\alpha_1^2+\alpha_3^2)+\sum_{i=0}^3X_i\cdot X_i$$
for $\lambda\not=1$. 
From the real version of (7) we have
$$
\split
&\quad d\pmatrix \lambda\alpha_1\\ \lambda\alpha_3\\ X_0\\ X_1\\ X_2\\ X_3\endpmatrix\\
&=-\,\pmatrix 0&-2\alpha_2& -\lambda X_1& \lambda X_0& \lambda X_3& -\lambda X_2\\
2\alpha_2 & 0 & -\lambda X_3& \lambda X_2& -\lambda X_1& \lambda X_0\\
\lambda^{-1} X_1& \lambda^{-1}X_3& 0 & -\Gamma_1 & -\Gamma_2-\alpha_2 & -\Gamma_3\\
-\lambda^{-1}X_0& -\lambda^{-1}X_2& \Gamma_1 & 0 & -\Gamma_3 & \Gamma_2-\alpha_2\\
-\lambda^{-1}X_3 & \lambda^{-1}X_1 & \Gamma_2+\alpha_2 & \Gamma_3 & 0 & -\Gamma_1\\
\lambda^{-1}X_2 & -\lambda^{-1}X_0& \Gamma_3& -\Gamma_2+\alpha_2& \Gamma_1&0
\endpmatrix 
\pmatrix \lambda\alpha_1\\ \lambda\alpha_3\\ X_0\\ X_1\\ X_2\\ X_3\endpmatrix\,\,.
\endsplit\tag8$$
The big matrix in (8) is the Levi-Civita connection form of the metric $g^{\CY}_{\lambda}$ on $\Z$. 
WE write this connection form as $\Gamma_{\lambda}$. 
Write
$$
\Omega_{\lambda}=\pmatrix {\Omega_{\lambda}}^{-2}_{-2}& 
{\Omega_{\lambda}}^{-2}_{-1}& {\Omega_{\lambda}}^{-2}_{0}& 
{\Omega_{\lambda}}^{-2}_{1}& {\Omega_{\lambda}}^{-2}_{2}&{\Omega_{\lambda}}^{-2}_{3}\\
{\Omega_{\lambda}}^{-1}_{-2}& {\Omega_{\lambda}}^{-1}_{-1}& {\Omega_{\lambda}}^{-1}_{0}& {\Omega_{\lambda}}^{-1}_{1}& {\Omega_{\lambda}}^{-1}_{2}& {\Omega_{\lambda}}^{-1}_{3}\\
{\Omega_{\lambda}}^{0}_{-2}& {\Omega_{\lambda}}^{0}_{-1}& {\Omega_{\lambda}}^{0}_{0}& {\Omega_{\lambda}}^{0}_{1}& {\Omega_{\lambda}}^{0}_{2}& {\Omega_{\lambda}}^{0}_{3}\\
{\Omega_{\lambda}}^{1}_{-2}& {\Omega_{\lambda}}^{1}_{-1}& {\Omega_{\lambda}}^{1}_{0}& {\Omega_{\lambda}}^{1}_{1}& {\Omega_{\lambda}}^{1}_{2}& {\Omega_{\lambda}}^{1}_{3}\\
{\Omega_{\lambda}}^{2}_{-2}& {\Omega_{\lambda}}^{2}_{-1}& {\Omega_{\lambda}}^{2}_{0}& {\Omega_{\lambda}}^{2}_{1}& {\Omega_{\lambda}}^{2}_{2}& {\Omega_{\lambda}}^{2}_{3}\\
{\Omega_{\lambda}}^{3}_{-2}& {\Omega_{\lambda}}^{3}_{-2}& {\Omega_{\lambda}}^{3}_{-2}& {\Omega_{\lambda}}^{3}_{-2}& {\Omega_{\lambda}}^{3}_{-2}& {\Omega_{\lambda}}^{3}_{-2}
\endpmatrix$$
for the curvature form of the metric $g^{\CY}_{\lambda}$ on $\Z$. 
The second structure equation 
$$d\Gamma_{\lambda}+\Gamma_{\lambda}\wedge\Gamma_{\lambda}=\Omega_{\lambda}$$
computes the curvature form $\Omega_{\lambda}$. From (8) and the second structure equation 
we have by direct computation the following expressions : 
$$
\split
& {\Omega_{\lambda}}^{-2}_{-2}={\Omega_{\lambda}}^{-1}_{-1}=0\,\,,\\
& {\Omega_{\lambda}}^{-1}_{-2}=4\alpha_3\wedge\alpha_1+2(X_2\wedge X_0+X_3\wedge X_1)\,\,,\\
&  {\Omega_{\lambda}}^{0}_{-2}=\lambda^{-1}(X_0\wedge\alpha_1+X_2\wedge\alpha_3)\,\,,\quad 
 {\Omega_{\lambda}}^{0}_{-1}=\lambda^{-1}(X_0\wedge\alpha_3-X_2\wedge\alpha_1)\,\,,\\
&  {\Omega_{\lambda}}^{1}_{-2}=\lambda^{-1}(X_1\wedge\alpha_1+X_3\wedge\alpha_3)\,\,,\quad 
 {\Omega_{\lambda}}^{1}_{-1}=\lambda^{-1}(X_1\wedge\alpha_3-X_3\wedge\alpha_1)\,\,,\\
&  {\Omega_{\lambda}}^{2}_{-2}=\lambda^{-1}(-X_0\wedge\alpha_3+X_2\wedge\alpha_1)\,\,,\quad
 {\Omega_{\lambda}}^{2}_{-1}=\lambda^{-1}(X_0\wedge\alpha_1+X_2\wedge\alpha_3)\,\,,\\
&  {\Omega_{\lambda}}^{3}_{-2}=\lambda^{-1}(-X_1\wedge\alpha_3+X_3\wedge\alpha_1)\,\,,\quad
 {\Omega_{\lambda}}^{3}_{-1}=\lambda^{-1}(X_1\wedge\alpha_1+X_3\wedge\alpha_3)\,\,,\\
&  {\Omega_{\lambda}}^{-2}_{0}=\lambda(\alpha_1\wedge X_0+\alpha_3\wedge X_2)\,\,,\quad
 {\Omega_{\lambda}}^{-1}_{0}=\lambda(\alpha_3\wedge X_0-\alpha_1\wedge X_2)\,\,,\\
&  {\Omega_{\lambda}}^{-2}_{1}=\lambda(\alpha_1\wedge X_1+\alpha_3\wedge X_3)\,\,,\quad
 {\Omega_{\lambda}}^{-1}_{1}=\lambda(\alpha_3\wedge X_1-\alpha_1\wedge X_3)\,\,,\\
&  {\Omega_{\lambda}}^{-2}_{2}=\lambda(-\alpha_3\wedge X_0+\alpha_1\wedge X_2)\,\,,\quad
 {\Omega_{\lambda}}^{-1}_{2}=\lambda(\alpha_1\wedge X_0+\alpha_3\wedge X_2)\,\,,\\
&  {\Omega_{\lambda}}^{-2}_{3}=\lambda(-\alpha_3\wedge X_1+\alpha_1\wedge X_3)\,\,,\quad
 {\Omega_{\lambda}}^{-1}_{3}=\lambda(\alpha_1\wedge X_1+\alpha_3\wedge X_3)\,\,,\\
&  {\Omega_{\lambda}}^{0}_{0}= {\Omega_{\lambda}}^{1}_{1}= {\Omega_{\lambda}}^{2}_{2}
= {\Omega_{\lambda}}^{3}_{3}=0\,\,,\\
&  {\Omega_{\lambda}}^{1}_{0}=\Omega^1_0+X_0\wedge X_1+X_2\wedge X_3
-(d\alpha_1-2\alpha_2\wedge\alpha_3)\,\,,\\
&  {\Omega_{\lambda}}^{2}_{0}=\Omega^2_0+2X_3\wedge X_1
-(d\alpha_2-2\alpha_3\wedge\alpha_1)+d\alpha_2\,\,,\\
&  {\Omega_{\lambda}}^{3}_{0}=\Omega^3_0-X_2\wedge X_1+X_0\wedge X_3
-(d\alpha_3-2\alpha_1\wedge\alpha_2)\,\,,\\
&  {\Omega_{\lambda}}^{2}_{1}=\Omega^2_1-X_3\wedge X_0+X_1\wedge X_2
+(d\alpha_3-2\alpha_1\wedge\alpha_2)\,\,,\\
&  {\Omega_{\lambda}}^{3}_{1}=\Omega^3_1+2X_2\wedge X_0
-(d\alpha_2-2\alpha_3\wedge\alpha_1)+d\alpha_2\,\,,\\
&  {\Omega_{\lambda}}^{3}_{2}=\Omega^3_2+X_2\wedge X_3+X_0\wedge X_1
+(d\alpha_1-2\alpha_2\wedge\alpha_3)\,\,.
\endsplit\tag9$$
We recall the fact that the curvature form of self-dual Einstein 4-manifolds decomposes as
$$
\split
& \Omega^{\mu}_0={\wt\Omega}^{\mu}_0+{\Omega'}^{\mu}_0
=4X_{\mu}\wedge X_0+{\Omega'}^{\mu}_0\,\,,\\
& \Omega^{\eta}_{\nu}={\wt\Omega}^{\eta}_{\nu}+{\Omega'}^{\eta}_{\nu}
=4X_{\eta}\wedge X_{\nu}+{\Omega'}^{\eta}_{\nu}\,\,,
\endsplit\tag10
$$
where $(\mu,\eta,\nu)$ is any cyclic permutation of $(1,2,3)$, $\wt\Omega$ is the curvature 
form of $S^4$ (with the standard metric with sectional curvature identically $4$) and $\Omega'$ 
is the $\pmatrix W_+&0\\0&0\endpmatrix$-part (which has no contribution to the Ricci tensor). 

From (9), (10) and the fact that $\Omega'$-part does not contribute to the Ricci tensor, 
we can compute the Ricci tensor by direct calculation using the definition of the Ricci tensor 
$$
R(e_i,e_j)=\sum_{k=1}^ng(R(e_i,e_k)e_k,e_j)=\sum_{k=1}^ng(\Omega^j_k(e_i,e_k)e_j,e_j)
$$
where $\{e_1,\dots,e_n\}$ is a system of orthonormal frame on an $n$-dimensional Riemannian 
manifold and $\Omega$ is the curvature form with respect to this orthonormal frame. 

In the following computation of the Ricci tensor, we write 
$$\{\lambda^{-1}\xi_{-2},\lambda^{-1}\xi_{-1},\xi_0,\xi_1,\xi_2,\xi_3\}$$
for the oriented orthonormal frame dual to the orthonormal coframe 
$$\{\lambda\alpha_1,\lambda\alpha_3,X_0,X_1,X_2,X_3\}\,\,.$$
The result of the computation of the Ricci tensor is as follows : 
$$
\split
& \Ric_{\lambda}(\lambda^{-1}\xi_{-2},\lambda^{-1}\xi_{-2})
=\Ric_{\lambda}(\lambda^{-1}\xi_{-1},\lambda^{-1}\xi_{-1})=\frac{4}{\lambda^2}+4\,\,,\\
& \Ric_{\lambda}(\lambda^{-1}\xi_{-2},\lambda^{-1}\xi_{-1})=0\,\,,\\
& \Ric_{\lambda}(\lambda^{-1}\xi_{-2},\xi_{0})=\Ric_{\lambda}(\lambda^{-1}\xi_{-2},\xi_{1})
=\Ric_{\lambda}(\lambda^{-1}\xi_{-2},\xi_{2})=\Ric_{\lambda}(\lambda^{-1}\xi_{-2},\xi_{3})\\
&=\Ric_{\lambda}(\lambda^{-1}\xi_{-1},\xi_{0})=\Ric_{\lambda}(\lambda^{-1}\xi_{-1},\xi_{1})
\Ric_{\lambda}(\lambda^{-1}\xi_{-1},\xi_{2})=\Ric_{\lambda}(\lambda^{-1}\xi_{-1},\xi_{3})\\
& =0\,\,,\\
& \Ric_{\lambda}(\xi_0,\xi_0)= \Ric_{\lambda}(\xi_1,\xi_1)= \Ric_{\lambda}(\xi_2,\xi_2)
= \Ric_{\lambda}(\xi_3,\xi_3)=\frac{2}{\lambda^2}+6\,\,,\\
& \Ric_{\lambda}(\xi_0,\xi_1)=\Ric_{\lambda}(\xi_0,\xi_2)=\Ric_{\lambda}(\xi_0,\xi_3)\\
&=\Ric_{\lambda}(\xi_1,\xi_2)=\Ric_{\lambda}(\xi_1,\xi_3)=\Ric_{\lambda}(\xi_2,\xi_3)\\
&=0\,\,.\endsplit$$

Summing up the above computation, we have :

\proclaim{Theorem 1.2 (Theorem 0.2 (2))} The Ricci tensor of the Chow-Yang metric $$g_{\lambda}^{\CY}=\lambda^2(\alpha_1^2+\alpha_3^2)+\sum_{i=0}^3X_i^2$$
is given by the formula
$$\split
\Ric(g^{\CY}_{\lambda})&=4(1+\lambda^{2})(\alpha_1^2+\alpha_3^2)
+2(\lambda^{-2}+3)(X_0^2+X_1^2+X_2^2+X_3^2)\\
&=2\lambda^{-2}(1+3\lambda^2)\,g^{\CY}_{\sqrt{\frac{2\lambda^2(1+\lambda^2)}{1+3\lambda^2}}}
\,\,.\endsplit
\tag11$$
Therefore, the Chow-Yang metric $g^{\CY}_{\lambda}$ has positive Ricci curvature 
and the Ricci tensor defines a Chow-Yang metric with different parameter, i.e., the Ricci map 
(the map which associate to each metric its Ricci tensor) is identified with a map on the $(\lambda^2,\rho)$ plane defined by
$$\Ric\,:\,(\lambda^2,\rho) \mapsto \biggl(\frac{2\lambda^2(1+\lambda^2)}{1+3\lambda^2},
2\lambda^{-2}(1+3\lambda^2)\biggr)\,\,.$$
\endproclaim

\beginsection 2. Ricci Flow
\par

In this section we study the Ricci flow on $\Z$ with initial metric taken from the family of the 
Chow-Yang metrics. We consider the two parameter family $\Cal F^{\CY}$ 
of scaled Chow-Yang metrics 
$\{\rho g^{\CY}_{\lambda}\}_{\rho>0,\lambda>0}$ where $\rho g_{\lambda}$ is expressed as, say, 
$$\rho g^{\CY}_{\lambda}=\rho\{\lambda^2(\alpha_1^2+\alpha_3^2)+\sum_{i=0}^3X_i^2\}\,\,.$$
Then the family $\Cal F^{\CY}$ is invariant under 
the sum with nonnegative coefficients. Moreover, Theorem 1.2 implies that the Ricci tensor 
of a Chow-Yang metric is again contained in the family $\Cal F^{\CY}$. 
Therefore, the Ricci flow with initial metric in $\Cal F^{\CY}$ stays in $\Cal F^{\CY}$ as long as 
the solution exists as Riemannian metrics on $\Z$. From Theorem 1.2 we have the explicit form 
of the Ricci flow with initial metric in the family $\Cal F^{\CY}$ : 

\proclaim{Theorem 2.1 (Theorem 0.3 (2))} (1) The Ricci flow equation $\p_tg=-2\Ric_g$ 
on the twistor space with initial metric in $\Cal F^{\CY}$ is given by the system of 
ordinary differential equations
$$
\left\{\aligned
& \frac{d}{dt}(\rho(t)\lambda(t)^2)=-8(1+\lambda^2(t))\,\,,\\
& \frac{d}{dt}\rho(t)=-4(\lambda(t)^{-2}+3)\,\,.\endaligned\right.
\tag12$$
The trajectory of the equation (12) in the $(\lambda^2,\rho)$-plane (this is identified with 
the Ricci flow trajectory) is given by
$$
\rho=(\text{\rm const.})\,\frac{\lambda^2}{|1-\lambda^2|^4}\,\,.\tag13
$$

(2) For any initial metric at time $t=0$ in the family $\Cal F^{\CY}$, the system of ordinary 
differential equations (12) has a solution defined on $(-\infty,T)$, i.e. the solution is extended 
for all negative reals and extinct at some finite time $T$ (i.e., as $t\to T$ the solution shrinks 
the space and become extinct at time $T$). The extinction time $T$ depends on the choice 
of the initial metric. 
\medskip

(3) Suppose that $\rho(0)=1$. 
If $\lambda(0)=1$, then the metric remains K\"ahler-Einstein ($\lambda(t)\equiv 1$) 
and the solution evolves just by homothety $\rho(t)=1-8t$ (in this case $T=\frac{1}{8}$). 

If $\lambda(0)<1$, then $\lim_{t\to-\infty}\lambda(t)=1$, $\lim_{t\to-\infty}\rho(t)=\infty$, 
$\lim_{t\to T}\lambda(t)=0$ and $\lim_{t\to T}\rho(t)=0$. 

If $\lambda(0)>1$, then 
$\lim_{t\to-\infty}\lambda(t)=1$, $\lim_{t\to-\infty}\rho(t)=\infty$, $\lim_{t\to T}\lambda(t)=\infty$, 
$\lim_{t\to T}\rho(t)=0$ and $\lim_{t\to T}\rho(t)\lambda^2(t)=0$. 

Suppose that $\lambda(0)\not=1$. Then, as $t$ becomes larger in the future direction, 
the deviation $|1-\lambda(t)|$ of the solution from being K\"ahler-Einstein becomes larger 
as well. As $t$ becomes larger in the past direction, then the solution becomes backward 
asymptotic to the solution in the case of 
$\lambda(0)=1$, i.e., the K\"ahler-Einstein metric is the asymptotic soliton of the 
Ricci flow under consideration. 
\medskip

(4) Suppose that $\lambda(0)<1$. Then the Gromov-Hausdorff limit of the Ricci flow solution 
as $t\to T$, scaled with the factor $\rho(t)^{-1}$,  is the original self-dual positive Einstein metric on 
$M$. 
\medskip

(5) Suppose that $\lambda(0)>1$. Then the Gromov-Haussdorff limit of the Ricci flow solution 
as $t\to T$, scaled with the factor $\rho(t)^{-1}$,  is the sub-Riemmanian metric defined on the 
horizontal distribution of the twistor space $\Cal Z$ which projects isometrically to 
the original self-dual positive Einstein metric on $M$. 
\endproclaim
\noindent
{\it Proof.} Putting $\mu(t)=\lambda^2(t)$ and assuming $\mu\not=1$, we get from (12) the 
system of ordinary differential equations
$$
\left\{\aligned
& \frac{d\rho}{dt}=-4\mu^{-1}-12\,\,,\\
& \frac{d(\rho\mu)}{dt}=-8-8\mu\,\,.
\endaligned\right.
$$
Eliminating $dt$ from the above two equations, we have
$$
d\log\rho=d\{\log\mu-4\log(1-\mu)\}\,\,.
$$
The solution of this equation represents a family of curves in the $(\mu,\rho)$-plane. 
Explicitly, we have (13) : 
$$
\rho=(\text{\rm const.})\,\frac{\mu}{|1-\mu|^4}\,\,.
$$
The meaning of (13) is this : Each Ricci flow solution of (12) is identified with the oriented curve 
in the $(\mu,\rho)$ ($\mu=\lambda^2$) plane, each of whose defining equation is given by (13) 
with a specified positive constant determined by the initial metric. 
The proof of Theorem 2.1 is now direct from (13). \qed

\proclaim{Remark 2.2 (comparison with canonical deformation metrics, cf. [Be] and [D-M])}
{\rm Let $g^{\can}_{\lambda}=\lambda^2g_{\FS}+g_M$ ($\lambda>0$) be the canonical 
deformation metrics on $\Cal Z$, where $g_{\FS}$ is the Fubini-Study metric of the fiber 
(curvature $\equiv 4$), 
$g_M$ is a quaternion K\"ahler metric with scalar curvature $S=48$ and the sum is 
defined with respect to the horizontal distribution corresponding to $g_M$. 
Then O'Neill's formula implies 
$$\Ric(g^{\can}_{\lambda})=4(1+\lambda^4)g_{FS}+4(3-\lambda^2)g_M\,\,\biggl[
=4(3-\lambda^2)\,
g^{\can}_{\sqrt{\frac{1+\lambda^4}{3-\lambda^2}}}\quad \text{\rm if $\lambda^2<3$}\biggr]\,\,.$$
We have two consequences from this formula. First, the Ricci tensor of $g^{\can}_{\lambda}$ 
is positive if and only if $\lambda^2<3$. 
Second, $g^{\can}_{\lambda}$ is Einstein if and only if 
$\lambda^2=1$ or $\lambda^2=\frac1{2}$. 
The case $\lambda^2=1$ corresponds to the submersion metric coming from 
$\Cal P\rightarrow \Cal Z$ which turns out to be K\"ahler-Einstein. Another Einstein metric (corresponding to the case $\lambda^2=\frac1{2}$) is non-K\"ahler. 
The Ricci flow equation for the canonical deformation metrics becomes
$$\left\{\aligned
& \frac{d\lambda^2}{dt}=-\frac{8}{\rho}(\lambda^2-1)(2\lambda^2-1)\\
& \frac{d\rho}{dt}=-8(3-\lambda^2)\,\,.\endaligned\right.
$$
For instance, if $\lambda>1$, then $\lambda$ decreases as $t$ increases (i.e., $\lambda$ 
decreases along the Ricci flow). Therefore the behavior of the Ricci flow in the canonical 
deformation is very different from Theorem 2.1. }\endproclaim

\beginsection 3. Uniformization of Self-Dual Positive Einstein 4-Manifolds
\par

Analyzing the ancient solutions in Theorem 2.1 we prove the following : 

\proclaim{Theorem 3.1} Let $(M,g)$ be a self-dual positive Einstein 
4-manifold and $g_{\lambda}^{\CY}$ be the Chow-Yang metric. Then we have the limit formula
$$\lim_{\lambda\to 0}|\nabla^{g_{\lambda}^{\CY}}\Rm(g_{\lambda}^{\CY})|_{g_{\lambda}^{\CY}}
=0\,\,.$$
\endproclaim
\noindent
{\it Proof}. We fix a small positive number $\delta$ and pick a sequence of points 
$\{(1+\delta,\rho_k)\}_{k=1}^{\infty}$ in the $(\mu,\rho)$ plane ($\mu=\lambda^2$) satisfying 
the property $\lim_{k\to\infty}\rho_k=\infty$. Correspondingly, we consider the Ricci flow 
solutions $\{g_k(t)\}_{k=1}^{\infty}$ in Theorem 2.1 with initial metric 
$g_k(0)=\rho_kg_{1+\delta}^{\CY}$. 
Let $T_k$ be the time when the solution $g_k(t)$ passes the line $\rho=1$, i.e., 
$g_k(T_k)=g_{\lambda_k}^{\CY}$ where $\{\lambda_k\}_{k=1}^{\infty}$ is a sequence of 
positive numbers satisfying the property $\lim_{k\to\infty}\lambda_k=\infty$. 
We note that $T_k\to\infty$ as $k\to\infty$. 
To see this, we mote that $\rho_k\to\infty$ as $k\to\infty$. 
We have from (12) that $\frac{dt}{d\rho}=-\frac14\frac{\lambda^2}{1+3\lambda^2}$ 
and therefore
$$
T_k=\frac14\int_1^{\rho_k}\frac{\lambda^2}{1+3\lambda^2}d\rho\,\,.
$$
As the trajectory of $g_k$ is given by the equation 
$\rho=c_k\frac{\lambda^2}{|1-\lambda^2|^4}$ in the $(\lambda,\rho)$ plane where $c_k$ is a 
positive constant satisfying the condition $\lim_{k\to\infty}c_k=\infty$, we see that 
as $\rho_k$ becomes large, the contribution 
from large $\lambda$ (which is $\approx \frac13$) in the above integral becomes 
more dominant. Therefore we have $T_k\to\infty$ as $k\to\infty$. 

Write $\Rm_k(x,t)$ for the Riemann curvature tensor of the solution metric $g_k(x,t)$. 
It follows from the computation in \S2 we have
$$
|\Rm_k(x,t)| \leq \frac{c_1}{T_k-t+1} \tag14
$$
whenever $0\leq t\leq T_k$ where the constant $c_1$ is chosen uniformly in $k$. 
We would like to estimate the asymptotic behavior of $\max_{x\in \Z}|\nabla\Rm_k(x,t)|$ 
when $t\to\infty$. For this purpose, we introduce
$$
F(x,t)=t\,|\nabla\Rm(x,t)|+\beta(t)\,|\Rm(x,t)|^2$$
where the function $\beta(t)$ is chosen in the form
$$
\beta(t)=1+\frac{c\,t}{T_k-t+1} \tag15
$$
with constant $c>0$ satisfying the condition
$$
1+c_1\,t\,|\Rm|-2\beta(t)<0\,\,.
$$
The strategy is to apply the parabolic maximum principle to the evolution inequality 
satisfied by the function $F(x,t)$ (see [Ba], [Sh1,2] and also [C-K, Chapt 7]). 
Following the notation in [C-K, p.227], we have
$$
\frac{\p}{\p t}|\Rm|^2=\triangle|\Rm|^2-2|\nabla \Rm|^2+(\Rm)^{*3}
$$
and
$$
\frac{\p}{\p t}|\nabla\Rm|^2=\triangle|\Rm|^2-2|\nabla^2\Rm|^2+\Rm*(\nabla\Rm)^{*2}\,\,.
$$
Combining the above two parabolic equations and (15), we conclude that 
the function $F(x,t)$ satisfies the parabolic inequality 
$$
\split
\frac{\p F}{\p t}&=|\nabla\Rm|^2+t\frac{\p}{\p t}|\nabla\Rm|^2
+\beta'(t)|\Rm|^2+\beta(t)\frac{\p}{\p t}|\Rm|^2\\
&=\triangle F+(|\nabla\Rm|^2+t\,\Rm*(\Rm)^{*2}-2\beta(t)|\nabla\Rm|^2)\\
&\quad +\beta(t)(\Rm)^{*3}+\beta'(t)|\Rm|^2-2t|\nabla^2\Rm|^2\\
& \leq \triangle F+(1+c_1t\,|\Rm|-2\beta(t))\,|\nabla\Rm|^2\\
& \quad +\beta(t)|\nabla\Rm|^3+\beta'(t)|\Rm|^2-2t|\nabla^2\Rm|^2\\
& \leq \triangle F+c_2\beta(t)|\Rm|^3+\beta'(t)|\Rm|^2
\endsplit\tag16
$$
where $c_2$ is a positive constant which can be chosen uniformly in $k$. 
Applying the parabolic maximum principle (see, for instance, [C-K, pp.95-96]) to (16) 
with the help of (14) and (15), we have
$$
\split
\sup_{x\in \Z}F(x,T_k) & \leq C\,\biggl\{
\beta(0)\biggl(\frac{c_1}{T_k+1}\biggr)^2+\int_0^{T_k}
\biggl(\frac{T_k}{(T_k-t+1)^3}+\frac{T_k}{(T_k-t+1)^4}\biggr)\,dt\biggr\}\\
& \leq \frac{C}{T_k}\,\,.
\endsplit
$$
Therefore we have
$$
|\nabla\Rm(T_k)|^2\leq \frac1{T_k}\sup_{x\in\Z}F(x,T_k)\leq \frac{C}{T_k^2}\to 0
$$
as $k\to\infty$. This implies
$$
\lim_{\lambda\to\infty}|\nabla^{g_{\lambda}^{\CY}}\Rm(g_{\lambda}^{\CY})|_{g_{\lambda}^{CY}} =0
$$
because $g_k=g_{\lambda_k}^{\CY}$ can be chosen arbitrary as far as 
$\lim_{k\to\infty}\lambda_k=\infty$. 
\qed
\medskip

We use the limit formula in Theorem 3.1 to prove Theorem 0.4 (Theorem 3.2). 
Intuitively, the behavior of $\nabla^{g_{\lambda}^{\CY}}\Rm(g_{\lambda}^{\CY})$ 
reduces that of $\nabla\Rm$ for the original $(M,g)$ because the fiber $\PP^1$ 
of the twistor fibration $\Z\rightarrow M$ blows up in the limit $\lambda\to\infty$. 

\proclaim{Theorem 3.2 (Theorem 0.4)} Any locally irreducible self-dual positive Einstein 
4-manifold is isometric to either $S^4$ with standard metric or 
$\PP^2(\C)$ with Fubini-Study metric. \endproclaim
\noindent
{\it Proof}. We choose an orthonormal basis $(e_A)$ of the tangent space 
$(T_mM,g_m)$ ($m\in M$) 
and extend it to an orthonormal frame on a neighborhood of $m$ by parallel transportation along  geodesics emanating from $m$. This defines a $4$-dimensional surface $S$ centered at 
$(e_A)\in \Cal P$ ($\Cal P$ being the holonomy reduction of the principal bundle of orthonormal 
frames of $M$) which is transversal to the vertical foliation. This determines a $4$-dimensional 
surface $S'$ in the twistor space $\Cal Z$ centered at a point $\wt m$ on a $\PP^1$-fiber over $m$, 
which is transversal to the $\PP^1$-fibration. The covariant derivative of the curvature 
tensor at $m$ is computed by differentiating the components of the curvature tensor w.r.to the 
orthonormal frames represented by points of $S$ (identified with $S'$ in $\Cal Z$) in the direction 
of a horizontal tangent vector at $\wt m$ (identified with a tangent vector of $M$ at $m$). 
In \S3, we computed the curvature form of the metric 
$g_{\lambda}=\lambda^2(\alpha_1^2+\alpha_3^2)+\sum_{i=0}^3{}^tX_i\cdot X_i$ 
on $\Cal Z$. For our purpose, 
we need:
$$\split
& {\Omega_{\lambda}}^0_0=\Omega^0_0-X_1\wedge X_1-X_3\wedge X_3,\quad 
{\Omega_{\lambda}}^1_1=\Omega^0_0-X_0\wedge X_0-X_2\wedge X_2\\
& {\Omega_{\lambda}}^2_2=\Omega^0_0-X_1\wedge X_1-X_3\wedge X_3,\quad 
{\Omega_{\lambda}}^3_3=\Omega^0_0-X_0\wedge X_0-X_2\wedge X_2\\
\endsplit$$
and
$$\split
& {\Omega_{\lambda}}^1_0=\Omega^1_0+X_0\wedge X_1+X_2\wedge X_3
-(d\alpha_1-2\alpha_2\wedge\alpha_3)\\
& {\Omega_{\lambda}}^2_0=\Omega^2_0+X_3\wedge X_1-X_1\wedge X_3
-(d\alpha_2-2\alpha_3\wedge\alpha_1)+d\alpha_2\\
& {\Omega_{\lambda}}^3_0=\Omega^3_0-X_2\wedge X_1+X_0\wedge X_3
-(d\alpha_3-2\alpha_1\wedge\alpha_2)\\
& {\Omega_{\lambda}}^2_1=\Omega^2_1-X_3\wedge X_0+X_1\wedge X_2
+(d\alpha_3-2\alpha_1\wedge\alpha_2)\\
& {\Omega_{\lambda}}^3_1=\Omega^3_1+X_2\wedge X_0-X_0\wedge X_2
-(d\alpha_2-2\alpha_3\wedge\alpha_1)+d\alpha_2\\
& {\Omega_{\lambda}}^3_2=\Omega^3_2+X_2\wedge X_3+X_0\wedge X_1
+(d\alpha_1-2\alpha_2\wedge\alpha_3)\,\,.
\endsplit$$
Taking the component in the $X_i$ ($i=0,1,2,3$) direction of the curvature tensor and 
taking the covariant derivative in the  $X_i$ ($i=0,1,2,3$) direction, we immediately 
conclude that the covariant derivatives of the  $X_i$ ($i=0,1,2,3$) part of the curvature tensor 
of the metric $g_{\lambda}$ of the twistor space $\Cal Z$ at $\wt m$ in the horizontal direction 
is equal to the covariant derivative in the corresponding direction 
of the curvature tensor of the quaternion K\"ahler manifold $(M,g)$ 
under question. On the other hand, we have from Theorem 4.2 the limit formula 
$$\lim_{\lambda\to\infty}|\nabla^{g^{\CY}_{\lambda}}
\Rm(g^{\CY}_{\lambda})|_{g^{\CY}_{\lambda}}=0\,\,.$$
This implies that the curvature tensor of the positive quaternion K\"ahler manifold 
$(M,g)$ must satisfy the condition $\nabla R\equiv 0$ 
from the beginning. This implies that $(M,g)$ is a symmetric space. Since we assumed that 
$(M,g)$ is irreducible, $(M,g)$ must be isometric to one of the locally irreducible 4-dimensional 
compact Riemannian symmetric spaces. \qed

\proclaim{Remark 3.3}{\rm Here we explain the ``intuition" behind Theorem 3.1 and 3.2. 
Let us recall that Perelman's No Local Collapsing Theorem [P] asserts 
the following. Suppose that a finite time singularity occurs in the space-time of the Ricci flow 
on a closed manifold. 
Then, as a limit of suitable parabolic scalings, we get a Ricci flow ancient solution which encodes 
all information of the singularity. From this view point, Theorem 2.1 means the following. 
Suppose that $\lambda>1$. Then the collapse of the twistor space where the base $M$ 
shrinks faster is realized as a finite time singularity arising in the Ricci flow starting at 
the metric $\rho g^{\CY}_{\lambda}$ ($\lambda>1$). 
Moreover, this Ricci flow coincides with the above mentioned ancient solution, which encodes all information of the collapse of $\Z$ where the base $M$ shrinks faster. 
Moreover, the shrinking family of K\"ahler-Einstein metrics on $\Z$ is its asymptotic 
soliton. Therefore, we can view the K\"ahler-Einstein metric of $\Z$ not only as a $t\to\infty$ 
limit (thermo-dynamical equilibrium) of the ``heat"-type equation (``stability" under the 
K\"ahler-Ricci flow of the positive K\"ahler-Einstein metric of Fano manifolds 
proved by Perelman and Tian-Zhu [T-Z]), but as a $t\to-\infty$ limit of the (non-K\"ahler) 
Ricci flow ancient solution. 
The $t\to-\infty$ limit should correspond to the ``minimum entropy". 
If we recall the ``intuitive" meaning of the thermo-dynamical entropy, we can expect that 
we can extract the local information of the K\"ahler-Einstein metric on $\Z$ (for instance, 
the ``local symmetry" of the original self-dual positive Einstein metrric on the base $M$).}
\endproclaim

\enddocument